\documentclass[12pt]{article}

\usepackage{setspace}

\usepackage{amsthm, amssymb, amstext}
\usepackage[fleqn]{amsmath}
\usepackage{latexsym}
\usepackage[dvips]{graphicx}
\usepackage{comment}
\usepackage{hyperref}
\usepackage{mathtools}
\usepackage{enumerate} 

\usepackage{todonotes}

\def\R{\mathcal{R}}

\newtheorem{theorem}{Theorem}

\theoremstyle{definition}

\usepackage{amsthm}

\title{Mixing colourings in $2K_2$-free graphs}
\author{Carl Feghali\thanks{Computer Science Institute of Charles University, Prague, Czech Republic, email: \texttt{feghali.carl@gmail.com} } \quad Owen Merkel\thanks{Department of Mathematics, Wilfrid Laurier University, Waterloo, Ontario, Canada.  Email: \texttt{owenmerkel@gmail.com}.}}

\date{}

\begin{document}
\maketitle

\begin{abstract}
The reconfiguration graph for the $k$-colourings of a graph $G$, denoted $\R_{k}(G)$, is the graph whose vertices are the $k$-colourings of $G$ and two colourings are joined by an edge if they differ in colour on exactly one vertex. For any $k$-colourable $P_4$-free graph $G$, Bonamy and Bousquet proved that $\R_{k+1}(G)$ is connected. In this short note, we complete the classification of the connectedness of $\R_{k+1}(G)$ for a $k$-colourable graph $G$ excluding a fixed path, by constructing  a $7$-chromatic $2K_2$-free (and hence $P_5$-free) graph admitting a frozen $8$-colouring. This settles a question of the second author. 
 \end{abstract}
 
\section{Introduction and result}
Let $G$ be a graph with vertex set $V(G)$ and edge set $E(G)$. For a positive integer $k$, a \emph{$k$-colouring} of $G$ is a mapping $\alpha \colon V(G) \to \{1, 2, \ldots, k\}$ such that $\alpha(u) \neq \alpha(v)$ whenever $uv \in E(G)$. The \emph{reconfiguration graph} for the $k$-colourings of $G$, denoted $\mathcal{R}_k(G)$, is the graph whose vertices are the $k$-colourings of $G$ and two colourings are joined by an edge if they differ in colour on exactly one vertex. We say that $G$ is \emph{$k$-mixing} if $\mathcal{R}_k(G)$ is connected.

A $k$-colouring $\varphi$ of $G$ is called \emph{frozen} if $\varphi$ is an isolated vertex in $R_k(G)$. In other words, $\varphi$ is frozen if for every vertex $v \in V(G)$, each of the $k$ colours appears in the closed neighbourhood of $v$. One technique to prove that a graph $G$ is not $k$-mixing is to exhibit a frozen $k$-colouring of $G$.

Let $k$ and $t$ be positive integers, and let $G$ be a $k$-colourable $P_t$-free graph. The subject of this note is to determine the smallest integer $\ell \geq k$ such that $G$ is $\ell$-mixing. Notice that if $G$ is the complete graph on $k$ vertices, $R_k(G)$ is disconnected, so we assume from now on that $\ell > k$. 

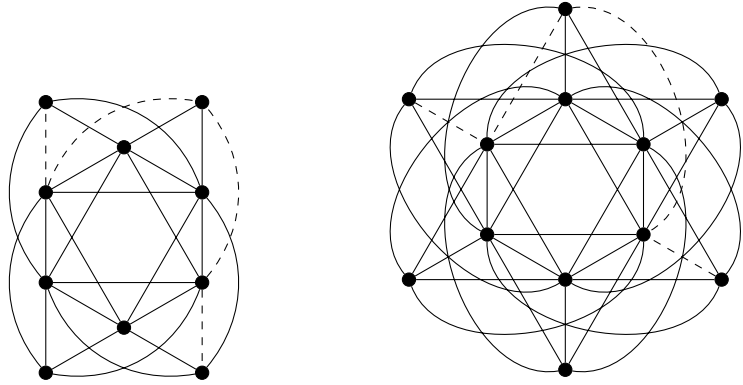
\begin{figure}
\centering
\begin{tikzpicture}[scale=1.2]
\tikzstyle{vertex}=[circle, draw, fill=black, inner sep=0pt, minimum size=5pt]
	
	\node[vertex](1) at (0,1) {};
    \node[vertex](2) at (cos{30},sin{30}) {};
    \node[vertex](3) at (cos{30},-sin{30}) {};
	\node[vertex](4) at (0,-1) {};
	\node[vertex](5) at (-cos{30},-sin{30}) {};
	\node[vertex](6) at (-cos{30},sin{30}) {};

	\node[vertex](7) at (cos{30},1.5) {};
	\node[vertex](8) at (cos{30},-1.5) {};
	\node[vertex](9) at (-cos{30},-1.5) {};
    \node[vertex](10) at (-cos{30},1.5) {};
    
    \draw 
    (1)--(2)
    (1)--(3)
    (1)--(5)
    (1)--(6)
    (1)--(7)
    (1)--(10)
    (2)--(3)
    (2)--(4)
    (2)--(6)
    (2)--(7)
    (2)to[bend left=40](8)
    (2)to[bend right=40](10)
    (3)--(4)
    (3)--(5)
    (3)to[bend left=40](9)
    (4)--(5)
    (4)--(6)
    (4)--(8)
    (4)--(9)
    (5)--(6)
    (5)to[bend right=40](8)
    (5)--(9)
    (5)to[bend left=40](10)
    (6)to[bend right=40](9);
    \draw[dashed]
    (6)--(10)
    (6)to[bend left=40](7)
    (3)to[bend right=40](7)
    (3) to (8);
    
\end{tikzpicture}
\hspace{10mm}
\begin{tikzpicture}[scale=1.2]
\tikzstyle{vertex}=[circle, draw, fill=black, inner sep=0pt, minimum size=5pt]

    \node[vertex](1) at (0,1) {};
    \node[vertex](2) at (cos{30},sin{30}) {};
    \node[vertex](3) at (cos{30},-sin{30}) {};
	\node[vertex](4) at (0,-1) {};
	\node[vertex](5) at (-cos{30},-sin{30}) {};
	\node[vertex](6) at (-cos{30},sin{30}) {};
	\node[vertex](7) at (0,2) {};
	\node[vertex](8) at (2*cos{30},2*sin{30}) {};
	\node[vertex](9) at (2*cos{30},-2*sin{30}) {};
	\node[vertex](10) at (0,-2) {};
	\node[vertex](11) at (-2*cos{30},-2*sin{30}) {};
	\node[vertex](12) at (-2*cos{30},2*sin{30}) {};
	
    \draw 
    (1) to (2)
    (1) to (3)
    (1) to (5)
    (1) to (6)
    (2) to (3)
    (2) to (4)
    (2) to (6)
    (3) to (4)
    (3) to (5)
    (4) to (5)
    (4) to (6)
    (5) to (6)
    (7) to (1)
    (7) to (2)
    (7) to[bend right=80] (5)
    (8) to (1)
    (8) to (2)
    (8) to (3)
    (8) to[bend right=80] (6)
    (8) to[bend left=80] (4)
    (9) to (2)
    (9) to (4)
    (9) to[bend right=80] (1)
    (9) to[bend left=80] (5)
    (10) to (3)
    (10) to (4)
    (10) to (5)
    (10) to[bend right=80] (2)
    (10) to[bend left=80] (6)
    (11) to (4)
    (11) to (5)
    (11) to (6)
    (11) to[bend right=80] (3)
    (11) to[bend left=80] (1)
    (12) to (5)
    (12) to (1)
    (12) to[bend right=80] (4)
    (12) to[bend left=80] (2);
    
    \draw[dashed]
    (12) to (6)
    (7) to (6)
    (7) to[bend left=80] (3)
    (9) to (3);
    
\end{tikzpicture}
\caption{Graphs mistaken to be $P_5$-free by Bonamy and Bousquet \cite{bonamy2018}. An induced $P_5$ is marked with dashed edges.}
\label{fig:notp5free}
\end{figure}

Bonamy and Bousquet \cite{bonamy2018} proved that any $k$-colourable $P_4$-free graph is $\ell$-mixing. Moreover,
using an example of Cereceda, van den Heuvel, and Johnson \cite{cereceda2008}, they also pointed out that, for every $\ell \geq 3$, there exists a $2$-colourable $P_6$-free graph with a frozen $\ell$-colouring. In that same paper, they also attempted to find a family of $k$-colourable $P_5$-free graphs with a frozen $\ell$-colouring. Unfortunately, as noted in \cite{merkel2020} the graphs they considered are not $P_5$-free (see Figure \ref{fig:notp5free}). This leaves open the case $t = 5$, which has since received some attention and has been considered for several subclasses of $P_5$-free graphs \cite{biedl2021, bonamy2018, feghali2020, merkel2020, merkel2021}. The aim of this short note is to settle this case, thereby answering a question raised in \cite{merkel2020}.

\begin{theorem}
\label{thm:main}
For every positive integer $p$, there exists a $7p$-colourable $P_5$-free graph that is not $8p$-mixing. 
\end{theorem}

We do this by constructing, for each positive integer $p$,  a $7p$-chromatic $2K_2$-free (and hence $P_5$-free) graph admitting a frozen $8p$-colouring as follows.

\begin{figure}
\centering
\begin{tikzpicture}[scale=2.75]
\tikzstyle{vertex}=[circle, draw, fill=black, inner sep=0pt, minimum size=5pt]
    
    \draw (0,0) circle (1);
    
    \node[vertex, label=above:1](0) at (0,1) {};
    \node[vertex, label=above:2](1) at (cos{67.5},sin{67.5}) {};
    \node[vertex, label=above:3](2) at (cos{45},sin{45}) {};
    \node[vertex, label=above:4](3) at (cos{22.5}, sin{22.5}) {};
    \node[vertex, label=right:5](4) at (1,0) {};
    \node[vertex, label=below:7](5) at (cos{22.5}, -sin{22.5}) {};
    \node[vertex, label=below:2](6) at (cos{45},-sin{45}) {};
    \node[vertex, label=below:3](7) at (cos{67.5},-sin{67.5}) {};
    \node[vertex, label=below:4](8) at (0,-1) {};
    \node[vertex, label=left:5](9) at (-cos{67.5},-sin{67.5}) {};
     \node[vertex, label=left:1](10) at (-cos{45},-sin{45}) {};
    \node[vertex, label=left:2](11) at (-cos{22.5}, -sin{22.5}) {};
    \node[vertex, label=left:3](12) at (-1,0) {};
    \node[vertex, label=above:4](13) at (-cos{22.5},sin{22.5}) {};
    \node[vertex, label=above:6](14) at (-cos{45},sin{45}) {};
    \node[vertex, label=above:7](15) at (-cos{67.5},sin{67.5}) {};
    
    \draw(0)--(2);
    \draw(0)--(3);
    \draw(0)--(4);
    \draw(0)--(5);
    \draw(0)--(6);
    \draw(0)--(7);
    \draw(0)--(9);
    \draw(0)--(11);
    \draw(0)--(12);
    \draw(0)--(13);
    \draw(0)--(14);
    \draw(1)--(3);
    \draw(1)--(4);
    \draw(1)--(5);
    \draw(1)--(7);
    \draw(1)--(8);
    \draw(1)--(10);
    \draw(1)--(12);
    \draw(1)--(13);
    \draw(1)--(14);
    \draw(1)--(15);
    \draw(2)--(4);
    \draw(2)--(5);
    \draw(2)--(6);
    \draw(2)--(8);
    \draw(2)--(9);
    \draw(2)--(11);
    \draw(2)--(13);
    \draw(2)--(14);
    \draw(2)--(15);
    \draw(3)--(5);
    \draw(3)--(6);
    \draw(3)--(7);
    \draw(3)--(9);
    \draw(3)--(10);
    \draw(3)--(12);
    \draw(3)--(14);
    \draw(3)--(15);
    \draw(4)--(6);
    \draw(4)--(7);
    \draw(4)--(8);
    \draw(4)--(10);
    \draw(4)--(11);
    \draw(4)--(13);
    \draw(4)--(14);
    \draw(4)--(15);
    \draw(5)--(7);
    \draw(5)--(8);
    \draw(5)--(9);
    \draw(5)--(10);
    \draw(5)--(11);
    \draw(5)--(12);
    \draw(5)--(14);
    \draw(6)--(8);
    \draw(6)--(9);
    \draw(6)--(10);
    \draw(6)--(12);
    \draw(6)--(13);
    \draw(6)--(15);
    \draw(7)--(9);
    \draw(7)--(10);
    \draw(7)--(11);
    \draw(7)--(13);
    \draw(7)--(14);
    \draw(8)--(10);
    \draw(8)--(11);
    \draw(8)--(12);
    \draw(8)--(14);
    \draw(8)--(15);
    \draw(9)--(11);
    \draw(9)--(12);
    \draw(9)--(13);
    \draw(9)--(14);
    \draw(9)--(15);
    \draw(10)--(12);
    \draw(10)--(13);
    \draw(10)--(14);
    \draw(10)--(15);
    \draw(11)--(13);
    \draw(11)--(14);
    \draw(11)--(15);
    \draw(12)--(14);
    \draw(12)--(15);
    \draw(13)--(15);

\end{tikzpicture}
\hspace{0mm}
\begin{tikzpicture}[scale=2.75]
\tikzstyle{vertex}=[circle, draw, fill=black, inner sep=0pt, minimum size=5pt]
    
    \draw (0,0) circle (1);
    
    \node[vertex, label=above:1](0) at (0,1) {};
    \node[vertex, label=above:2](1) at (cos{67.5},sin{67.5}) {};
    \node[vertex, label=above:3](2) at (cos{45},sin{45}) {};
    \node[vertex, label=above:4](3) at (cos{22.5}, sin{22.5}) {};
    \node[vertex, label=right:5](4) at (1,0) {};
    \node[vertex, label=below:6](5) at (cos{22.5}, -sin{22.5}) {};
    \node[vertex, label=below:7](6) at (cos{45},-sin{45}) {};
    \node[vertex, label=below:8](7) at (cos{67.5},-sin{67.5}) {};
    \node[vertex, label=below:1](8) at (0,-1) {};
    \node[vertex, label=left:2](9) at (-cos{67.5},-sin{67.5}) {};
     \node[vertex, label=left:3](10) at (-cos{45},-sin{45}) {};
    \node[vertex, label=left:4](11) at (-cos{22.5}, -sin{22.5}) {};
    \node[vertex, label=left:5](12) at (-1,0) {};
    \node[vertex, label=above:6](13) at (-cos{22.5},sin{22.5}) {};
    \node[vertex, label=above:7](14) at (-cos{45},sin{45}) {};
    \node[vertex, label=above:8](15) at (-cos{67.5},sin{67.5}) {};
    
    \draw(0)--(2);
    \draw(0)--(3);
    \draw(0)--(4);
    \draw(0)--(5);
    \draw(0)--(6);
    \draw(0)--(7);
    \draw(0)--(9);
    \draw(0)--(11);
    \draw(0)--(12);
    \draw(0)--(13);
    \draw(0)--(14);
    \draw(1)--(3);
    \draw(1)--(4);
    \draw(1)--(5);
    \draw(1)--(7);
    \draw(1)--(8);
    \draw(1)--(10);
    \draw(1)--(12);
    \draw(1)--(13);
    \draw(1)--(14);
    \draw(1)--(15);
    \draw(2)--(4);
    \draw(2)--(5);
    \draw(2)--(6);
    \draw(2)--(8);
    \draw(2)--(9);
    \draw(2)--(11);
    \draw(2)--(13);
    \draw(2)--(14);
    \draw(2)--(15);
    \draw(3)--(5);
    \draw(3)--(6);
    \draw(3)--(7);
    \draw(3)--(9);
    \draw(3)--(10);
    \draw(3)--(12);
    \draw(3)--(14);
    \draw(3)--(15);
    \draw(4)--(6);
    \draw(4)--(7);
    \draw(4)--(8);
    \draw(4)--(10);
    \draw(4)--(11);
    \draw(4)--(13);
    \draw(4)--(14);
    \draw(4)--(15);
    \draw(5)--(7);
    \draw(5)--(8);
    \draw(5)--(9);
    \draw(5)--(10);
    \draw(5)--(11);
    \draw(5)--(12);
    \draw(5)--(14);
    \draw(6)--(8);
    \draw(6)--(9);
    \draw(6)--(10);
    \draw(6)--(12);
    \draw(6)--(13);
    \draw(6)--(15);
    \draw(7)--(9);
    \draw(7)--(10);
    \draw(7)--(11);
    \draw(7)--(13);
    \draw(7)--(14);
    \draw(8)--(10);
    \draw(8)--(11);
    \draw(8)--(12);
    \draw(8)--(14);
    \draw(8)--(15);
    \draw(9)--(11);
    \draw(9)--(12);
    \draw(9)--(13);
    \draw(9)--(14);
    \draw(9)--(15);
    \draw(10)--(12);
    \draw(10)--(13);
    \draw(10)--(14);
    \draw(10)--(15);
    \draw(11)--(13);
    \draw(11)--(14);
    \draw(11)--(15);
    \draw(12)--(14);
    \draw(12)--(15);
    \draw(13)--(15);
    
\end{tikzpicture}
\caption{A 7-colouring and a frozen 8-colouring of $G$.}
\label{fig:G}
\end{figure}
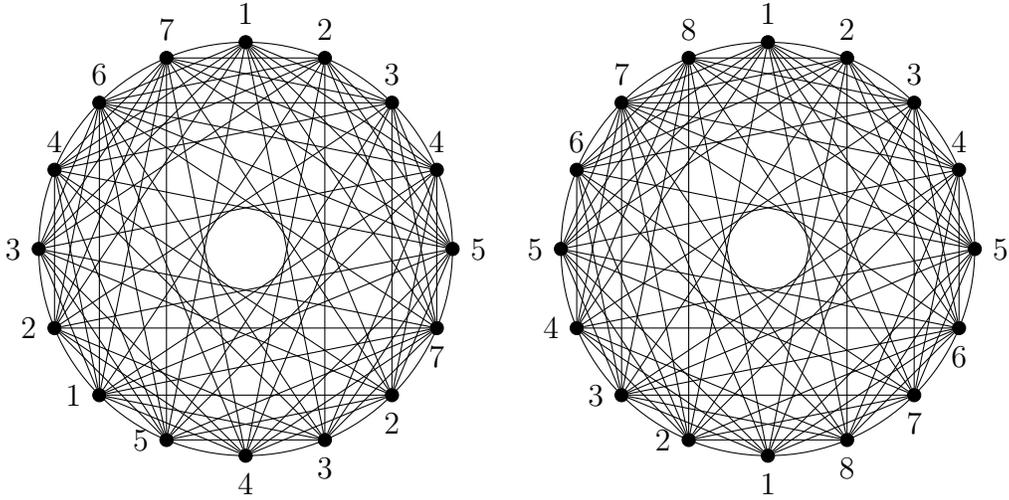

Let $C$ be a cycle on $16$ vertices and $H$ be the fourth power of $C$. Let $v$ be a vertex in $V(C) = V(H)$, and $\alpha$ and $\beta$ be, respectively, the colourings  $1234572345123467$ and $1234567812345678$ of  $V(H)$ in clockwise order around $C$ starting at $v$. Let $G$ be the graph obtained from $H$ by adding any missing edges between pairs of vertices that are coloured distinctly under \emph{both} $\alpha$ and $\beta$; see Figure \ref{fig:G}. By construction, $\alpha$ and $\beta$ are proper in $G$, and as every colour appears in the closed neighbourhood of every vertex in $H$ (and hence in $G$) under $\beta$, the colouring $\beta$ is frozen. It can be checked via a standard implementation in, for example, SageMath, that $G$ is $2K_2$-free. This establishes Theorem \ref{thm:main} when $p = 1$.  For $p \geq 2$, we take the pairwise complete join of $p$ copies of $G$ as our graph.  

\section{Closing remarks}

It has been shown that a $k$-colourable graph $G$ is $(k+1)$-mixing if $G$ is in one of several subclasses of $P_5$-free graphs. This includes when $k=2$ \cite{bonamy2014}, for co--chordal graphs, for ($P_5$, $\overline{P_5}$, $C_5$)-free graphs when $k=3$ \cite{feghali2020}, for $P_4$-sparse graphs \cite{biedl2021}, and for $3K_1$-free graphs \cite{merkel2021}. We leave as an open problem whether a $k$-colourable $P_5$-free graph is $(k+1)$-mixing for $3 \le k \le 6$.

In every known example of a $k$-colourable (arbitrary) graph $G$ that is not $(k+1)$-mixing, $G$ has a frozen $(k+1)$-colouring. It would be interesting to know if there are $k$-colourable graphs that are not $(k+1)$-mixing, but where every connected component of $\R_{k+1}(G)$ has more than one vertex.

\section*{Acknowledgments}
Carl Feghali was supported by grant 19-21082S of the Czech Science Foundation. Owen Merkel was partially supported by Natural Sciences and Engineering Research Council of Canada (NSERC) grant RGPIN-2016-06517.

\end{document}